\documentclass{amsart}
\usepackage{amsfonts}

\newtheorem{theorem}{Theorem}
\theoremstyle{plain}

\newtheorem{example}{Example}

\newtheorem{lemma}{Lemma}

\newtheorem{proposition}{Proposition}

\numberwithin{equation}{section}

\begin{document}
\title{ON THE CLASSIFICATION OF VECTOR BUNDLES WITH PERIODIC MAPS}
\author{ABDELOUAHAB AROUCHE}
\address{USTHB MATH. BP 32 EL ALIA 16111 BAB EZZOUAR ALGER ALGERIE.}
\email{abdarouche@hotmail.com.}
\date{March 10, 2005.}
\subjclass[2000]{ 19L47, 55R91.}
\keywords{Stiefel manifolds, equivariant bundles, completion}
\maketitle

\begin{abstract}
We give an explicit decription for univeral principal U(r)-bundles with
periodic map by means of equivariant Stiefel manifolds. We then show that
the associated equivariant vector bundle is equivalent to the canonical one
given by G. Segal. Finally, we investigate some ideals involved in the
equivariant K-theory of this classifying space.
\end{abstract}

\section{Introduction}

For a compact Lie group $\Gamma $ and a topological group $G$, we call a $%
(\Gamma ;G)$-bundle any principal $G$-bundle on which $\Gamma $ acts by
bundle maps [D], [L]. It has been shown in [D] and [L] that a universal $%
(\Gamma ;G)$-bundle is given by $E(\Gamma ;G)=\underset{H\in F}{*}\frac{%
\Gamma \times G}{H}$, where $F=F\left( \Gamma ;G\right) $ is the family of
subgroups $H\leq \Gamma \times G$ such that $H\cap G=1$ (see also [AHJM]).
This space is characte rized by the fact that the fixed points set $E(\Gamma
;G)^{H}$ is contractible if $H\in F$ and empty otherwise [J].

\smallskip

\section{Two lemmas}

In the sequel, we shall need the following

\begin{lemma}
if $G=G_{1}\times G_{2}$, then there is a $(\Gamma \times G)$-homotopy
equivalence :
\begin{equation*}
E(\Gamma ;G)\simeq E(\Gamma ;G_{1})\times E(\Gamma ;G_{2}).
\end{equation*}
\end{lemma}

\begin{lemma}
If $\Gamma _{1}\leq \Gamma _{2}$ and $G_{1}\leq G_{2}$, then there exists a $%
(\Gamma _{1}\times G_{1})$-homotopy equivalence :
\begin{equation*}
E(\Gamma _{1};G_{1})\simeq res_{(\Gamma _{1}\times G_{1})}E(\Gamma
_{2};G_{2}).
\end{equation*}
\end{lemma}

\smallskip

\section{The main theorem}

In what follows, we restrict our study to the case $\Gamma =\mathbf{Z}_{n}$,
the cyclic group of order $n$. Recall that the Stiefel manifolds are defined
by
\begin{equation*}
V_{r,s}=\left\{ (v_{1},...,v_{r}),v_{i}\in \mathbf{C}^{s},v_{i}.v_{j}=\delta
_{ij},i,j=1,...,r\right\}
\end{equation*}

and
\begin{equation*}
V_{r,\infty }=\underset{\overset{s}{\rightarrow }}{\lim }V_{r,s}.
\end{equation*}
Let $\mathbf{Z}_{n}\times U(r)$ act on the Stiefel manifold $V_{r,\infty }$
by :

\begin{equation*}
(\gamma ,a).(v_{1},...,v_{r})=(w_{1},...,w_{r})
\end{equation*}

such that
\begin{equation*}
w_{j}^{k}=\gamma ^{k}.\sum_{i=1}^{r}\overline{a_{ji}}v_{i}^{k}\text{, }%
j=1,...,r\text{, }k=1,...,\infty
\end{equation*}

where $\gamma $ is a generator of $\mathbf{Z}_{n}$, $a\in U(r)$, and $%
(v_{1},...,v_{r})\in V_{r,\infty }$.

We then have :

\begin{theorem}
There is a $\left( \mathbf{Z}_{n}\times U\left( r\right) \right) $-homotopy
equivalence :
\begin{equation*}
E\left( \mathbf{Z}_{n};U\left( r\right) \right) \simeq V_{r,\infty }
\end{equation*}
\end{theorem}

\begin{proof}
It is enough to show that $V_{r,\infty }$ endowed with the so defined $%
\left( \mathbf{Z}_{n}\times U\left( r\right) \right) $-action satisfies the
fixed points set characterization of $E\left( \mathbf{Z}_{n};U\left(
r\right) \right) $. Now, up to conjugation, we have :
\begin{equation*}
F\left( \mathbf{Z}_{n};U\left( r\right) \right) =\left\{ H_{d,s}:d\mid
n;s=(s_{1},...,s_{d});0\leq s_{i}\leq
r,i=1,...,d;\sum_{i=1}^{d}s_{i}=r\right\}
\end{equation*}%
with :
\begin{equation*}
H_{d,s}=\langle \lambda ,\rho \left( \lambda \right) =\left(
\begin{array}{cccccccccc}
\lambda &  &  &  &  &  &  &  &  &  \\
& \ddots &  &  &  &  &  &  &  &  \\
&  & \lambda &  &  &  &  &  &  &  \\
&  &  & \lambda ^{2} &  &  &  &  &  &  \\
&  &  &  & \ddots &  &  &  &  &  \\
&  &  &  &  & \lambda ^{2} &  &  &  &  \\
&  &  &  &  &  & \ddots &  &  &  \\
&  &  &  &  &  &  & \lambda ^{d}=1 &  &  \\
&  &  &  &  &  &  &  & \ddots &  \\
&  &  &  &  &  &  &  &  & \lambda ^{d}=1%
\end{array}%
\right) \rangle .
\end{equation*}%
The matrix $\rho \left( \lambda \right) $ contain $s_{i}$ diagonal
coefficients equal to $\lambda ^{i}$ ($\lambda $ is a primitive $d$-th root
of unity), all the other coefficient being zero. It i obvious that if%
\begin{equation*}
\left( 1,a\right) .\left( v_{1},...,v_{r}\right) =\left(
v_{1},...,v_{r}\right) ,
\end{equation*}%
then $a=\mathbf{I}_{r}$. Hence, $E\left( \mathbf{Z}_{n};U\left( r\right)
\right) ^{H}=\emptyset $ when $H\notin F\left( \mathbf{Z}_{n};U\left(
r\right) \right) $, on one hand. On the other hand, $\left(
v_{1},...,v_{r}\right) \in E\left( \mathbf{Z}_{n};U\left( r\right) \right)
^{H_{d,s}}$ iff
\begin{equation*}
\forall j=1,...,r:\forall k=1,...,\infty :v_{j}^{k}=\lambda
^{k-m_{j}}.v_{j}^{k};m_{j}=\rho \left( \lambda \right) _{jj}.
\end{equation*}%
Hence, $\forall j=1,...,r:v_{j}^{k}=0$ unless $k=m_{j} mod d]$. It i
not difficult to ssee that $E\left( \mathbf{Z}_{n};U\left( r\right)
\right) ^{H_{d,s}}$ is contractible (see $sec.6$. for an explicit
example).

\begin{example}
Consider the case $r=1$. We have $E\left( \mathbf{Z}_{n};\mathbf{S}%
^{1}\right) \simeq \mathbf{S}^{\infty }$. Using lemma 2.1, we get :
\begin{equation*}
E\left( \mathbf{Z}_{n};\mathbf{T}^{r}\right) \simeq \prod_{i=1}^{r}\mathbf{S}%
^{\infty }.
\end{equation*}
If we suppose that $G$ is a compact Lie group, we can embed it in a $U(r)$,
for some suitable $r$. Then, by use of lemma 2.2, there is a $\left( \mathbf{%
Z}_{n}\times G\right) $-action on $V_{r,\infty }$ and a $\left( \mathbf{Z}%
_{n}\times G\right) $-homotopy equivalence :
\begin{equation*}
E\left( \mathbf{Z}_{n};G\right) \simeq V_{r,\infty }.
\end{equation*}
\end{example}
\end{proof}

\section{The associated \textbf{Z}$_{n}$-vector bundle.}

Let $\mathbf{M}$ be the $\mathbf{Z}_{n}$-module defined by : $\mathbf{M=C}%
^{\infty }$ and%
\begin{equation*}
\gamma .\left( z_{1},z_{2},....\right) =\left( \gamma .z_{1},\gamma
^{2}.z_{2},....\right)
\end{equation*}%
$.$ Let $G_{r,\infty }$ be the Grassmannian manifold of $r$-dimensional
subspaces of $\mathbf{M}$. Then the right $U(r)$-action on $V_{r,\infty }$
defined by
\begin{equation*}
\left( v_{1},...,v_{r}\right) .a=\left( 1,a^{-1}\right) .\left(
v_{1},...,v_{r}\right)
\end{equation*}
has as a quotient space $V_{r,\infty }/ U(r)=G_{r,\infty }$. There
is a canonical $\mathbf{Z}_{n}$-vector bundle on $G_{r,\infty }$
whose totalspace is
\begin{equation*}
E_{\mathbf{M}}=\left\{ \left( N,z\right) :z\in N\right\} \subseteq
G_{r,\infty }\times \mathbf{M}
\end{equation*}
and which is universal [S].

\begin{theorem}
The $\mathbf{Z}_{n}$-vector bundle associated with the principal $\left(
\mathbf{Z}_{n};U(r)\right) $-bundle $V_{r,\infty }\rightarrow G_{r,\infty }$
is equivalent to $E_{\mathbf{M}}\rightarrow G_{r,\infty }$
\end{theorem}

\begin{proof}
The total space of the asociated $\mathbf{Z}_{n}$-vector bundle is $%
V_{r,\infty }\times _{U(r)}\mathbf{C}^{r}$, with the $\mathbf{Z}_{n}$-action
:%
\begin{equation*}
\gamma .\left[ \left( v_{1},...,v_{r}\right) ,y\right] =\left[ \left( \gamma
,\mathbf{I}_{r}\right) .\left( v_{1},...,v_{r}\right) ,y\right] .
\end{equation*}
Following [H, ch. 7, (7.1)], we define a map
\begin{equation*}
f:V_{r,\infty }\times _{U(r)}\mathbf{C}^{r}\rightarrow E_{\mathbf{M}},
\end{equation*}
by
\begin{equation*}
f\left[ \left( v_{1},...,v_{r}\right) ,y\right] =\left( \langle
v_{1},...,v_{r}\rangle ,\sum_{i=1}^{r}y_{i}v_{i}\right)
\end{equation*}%
. According to [H], $f$ is a vector bundle isomorphism. It is easy to check
its $\mathbf{Z}_{n}$-equivariance.
\end{proof}

\section{Characteristic classes}

Let $\Gamma $ be a compact Lie group and $G$ a topological group. Then $%
E\left( \Gamma ;G\right) $ is a right $G$-space by $e.g=\left( 1,g\right) .e$%
; we denote by $B\left( \Gamma ;G\right) $ the quotient space $E\left(
\Gamma ;G\right) /G$. According to [AHJM], the equivariant K-theory ring $%
K_{\Gamma }^{\ast }\left( B\left( \Gamma ;G\right) \right) $ is isomorphic
to the completion $R\left( \Gamma \times G\right) _{F}^{\symbol{94}}$ of the
complex representation ring with respect to the $F$-adic topology. Recall
that the $F$-adic topology is defined on $R\left( \Gamma \times G\right) $
by the set of ideals $I_{H}=\ker \left\{ R\left( \Gamma \times G\right)
\rightarrow R\left( H\right) \right\} $, for $H\in F$. In the case of $%
\Gamma =\mathbf{Z}_{2}$ and $G=\mathbf{T}^{r}$, we have :
\begin{equation*}
F=H_{k}:k=-1,0,...,r,
\end{equation*}%
where%
\begin{equation*}
H_{-1}=1,
\end{equation*}
\begin{equation*}
H_{0}=\left\{ 1,\left( -1,\mathbf{I}_{r}\right) \right\} ,H_{r}=\left\{
1,\left( -1,-\mathbf{I}_{r}\right) \right\}
\end{equation*}
and for $k=1,...,r-1:$%
\begin{equation*}
H_{k}=\left\{ 1,\left( -1,\left(
\begin{array}{ll}
-\mathbf{I}_{k} & 0 \\
0 & \mathbf{I}_{r-k}%
\end{array}%
\right) \right) \right\} .
\end{equation*}%
Putting%
\begin{equation*}
I_{k}=\ker \left\{ R\left( \mathbf{Z}_{2}\times \mathbf{T}^{r}\right)
\rightarrow R\left( H_{k}\right) \right\} ,
\end{equation*}
the $F$-adic topology is defined on $R\left( \mathbf{Z}_{2}\times \mathbf{T}%
^{r}\right) $ by the ideal :
\begin{equation*}
I=\bigcap_{k=-1}^{r}I_{k}.
\end{equation*}

\begin{proposition}
The the $F$-adic topology is defined on $R\left( \mathbf{Z}_{2}\times
\mathbf{T}^{r}\right) $ by the ideal :
\begin{equation*}
I=I_{0}\cong R\left( \mathbf{Z}_{2}\right) \otimes I_{\mathbf{T}^{r}},
\end{equation*}
where $I_{\mathbf{T}^{r}}$ denotes the augmentation ideal, and that on $%
R\left( \mathbf{Z}_{2}\times \mathbf{U}\left( r\right) \right) $ is defined
by :
\begin{equation*}
J=I\cap R\left( \mathbf{Z}_{2}\times \mathbf{U}\left( r\right) \right) .
\end{equation*}
In general, the $F$-adic topology on $R\left( \mathbf{Z}_{n}\times \mathbf{U}%
\left( r\right) \right) $ is defined by the ideals $I_{d,j}$ coming from the
subgroups :
\begin{equation*}
K_{d,j}=\langle \lambda ,\left(
\begin{array}{llll}
\lambda ^{j} &  &  &  \\
& 1 &  &  \\
&  & \ddots &  \\
&  &  & 1%
\end{array}
\right) \rangle ;\text{ }d\mid n;\text{ }j=1,...,d;\text{ }\left( \lambda
^{d}=1\right) .
\end{equation*}
\end{proposition}

\begin{proof}
Easy. The general case follows from the description of $F\left( \mathbf{Z}%
_{2};U\left( r\right) \right) $ given in thm (3.1).

\begin{example}
$K_{d,d}\cong \mathbf{Z}_{d}$, $I_{d,d}\cong R\left( \mathbf{Z}_{d}\right)
\otimes I_{U\left( r\right) }$ and $I_{d,d-1}\cong I_{\left( \mathbf{Z}%
_{d}\times \mathbf{T}^{r}\right) }\cap R\left( \mathbf{Z}_{d}\times U\left(
r\right) \right) $.
\end{example}
\end{proof}

\section{Homotopies in Stiefel manifolds}

Let $\mathbf{C}^{\infty }$ be the telescope $\bigcup_{n\geq 1}\mathbf{C}%
^{n}, $ and let
\begin{eqnarray*}
\left( \mathbf{C}^{\infty }\right) ^{even} &=&\left\{ v\in \mathbf{C}%
^{\infty }:v_{2n+1}=0;\quad \forall n=1,...,\infty .\right\} , \\
\left( \mathbf{C}^{\infty }\right) ^{odd} &=&\left\{ v\in \mathbf{C}^{\infty
}:v_{2n}=0;\quad \forall n=1,...,\infty .\right\} .
\end{eqnarray*}%
Two maps :
\begin{equation*}
g^{even},g^{odd}:\mathbf{C}^{n}\times \mathbf{I\rightarrow C}^{2n}
\end{equation*}%
are defined in the following way :
\begin{equation*}
g_{t}^{even}\left( v_{1},...,v_{n}\right) =\left( 1-t\right) \left(
v_{1},...,v_{n}\right) +t\left( 0,v_{1},...,0,v_{n}\right)
\end{equation*}%
and
\begin{equation*}
g_{t}^{odd}\left( v_{1},...,v_{n}\right) =\left( 1-t\right) \left(
v_{1},...,v_{n}\right) +t\left( v_{1},0,...,v_{n},0\right) .
\end{equation*}%
These maps extend to :
\begin{equation*}
g^{even},g^{odd}:\mathbf{C}^{\infty }\times \mathbf{I\rightarrow C}^{\infty
},
\end{equation*}%
and satisfy :
\begin{equation*}
g_{0}^{even}=g_{0}^{odd}=Id_{\mathbf{C}^{\infty }};
\end{equation*}%
\begin{equation*}
g_{1}^{even}\left( \mathbf{C}^{\infty }\right) \subseteq \left( \mathbf{C}%
^{\infty }\right) ^{even};g_{1}^{odd}\left( \mathbf{C}^{\infty }\right)
\subseteq \left( \mathbf{C}^{\infty }\right) ^{odd}.
\end{equation*}%
[H, ch. 3, (6.1)].

Similarly, we define
\begin{equation*}
\left( \mathbf{V}_{r,\infty }\right) ^{even}=\left\{ \left(
v_{1},...,v_{r}\right) \in \mathbf{V}_{r,\infty }:v_{j}^{2n+1}=0,\quad
\forall j=1,...,r;\quad \forall n=1,...,\infty \right\} ,
\end{equation*}
and
\begin{equation*}
\left( \mathbf{V}_{r,\infty }\right) ^{odd}=\left\{ \left(
v_{1},...,v_{r}\right) \in \mathbf{V}_{r,\infty }:v_{j}^{2n}=0,\quad \forall
j=1,...,r;\quad \forall n=1,...,\infty \right\} .
\end{equation*}

\begin{proposition}
Ther are maps :
\begin{equation*}
G^{even},G^{odd}:V_{r,\infty }\times \mathbf{I\rightarrow V}_{r,\infty },
\end{equation*}
such that
\begin{equation*}
G_{0}^{even}=G_{0}^{odd}=Id_{V_{r,\infty }};
\end{equation*}
\begin{equation*}
G_{1}^{even}\left( \mathbf{V}_{r,\infty }\right) \subseteq \left( \mathbf{V}%
_{r,\infty }\right) ^{even};G_{1}^{odd}\left( \mathbf{V}_{r,\infty }\right)
\subseteq \left( \mathbf{V}_{r,\infty }\right) ^{odd}.
\end{equation*}
\end{proposition}

\begin{proof}
we start by the case $r=1$. It is easy to see that $g_{t}^{even}(v)$ (resp. $%
g_{t}^{odd}\left( v\right) $) is not zero if $v$ is not. We define the
desired maps :
\begin{equation*}
G^{even},G^{odd}:\mathbf{S}^{\infty }\times \mathbf{I\rightarrow S}^{\infty
}.
\end{equation*}
by :
\begin{equation*}
G_{t}^{*}\left( v\right) =\frac{g_{t}^{*}\left( v\right) }{\left\|
g_{t}^{*}\left( v\right) \right\| },
\end{equation*}
for $*=even,odd$. In general, let $v\in \mathbf{C}^{n}$, and denote by $w$
the firt projection of $g_{t}^{*}\left( v\right) $ on $\mathbf{C}^{n}$. The
corresponding matrix is :
\begin{equation*}
\left(
\begin{array}{ccccccccccc}
1 &  &  &  &  &  &  &  &  &  &  \\
& 1-t &  &  &  &  &  &  &  &  &  \\
& t & 1-t &  &  &  &  &  &  &  &  \\
&  & 0 & 1-t &  &  &  &  &  &  &  \\
&  & t & 0 & 1-t &  &  &  &  &  &  \\
&  &  & 0 & 0 & 1-t &  &  &  &  &  \\
&  &  & t & 0 & 0 & 1-t &  &  &  &  \\
&  &  &  &  &  &  & \ddots &  &  &  \\
&  &  &  &  &  &  &  & 1-t &  &  \\
&  &  &  &  & t & 0 & \cdots & 0 & 1-t &  \\
&  &  &  &  &  &  &  &  &  & 1-t%
\end{array}
\right) .
\end{equation*}
We may delete the first or the last row and column according to the parity
of $n$ and $*$. We can prove by induction that this matrix determinant is $%
\left( 1-t\right) ^{n}$. We define the desired map by use of the
Gram-Schmidt map $GS$ :
\begin{equation*}
G_{t}^{*}\left( v_{1},...,v_{r}\right) =GS\left( g_{t}^{*}\left(
v_{1}\right) ,...,g_{t}^{*}\left( v_{r}\right) \right) .
\end{equation*}
\end{proof}

\bigskip

\bigskip

\bigskip

\end{document}